
\input amstex
\magnification=1200
\documentstyle{amsppt}
\NoRunningHeads
\NoBlackBoxes
\define\Mat{\operatorname{Mat}}
\define\so{\operatorname{\frak s\frak o}}
\define\Hom{\operatorname{Hom}}
\document\qquad\qquad\qquad\qquad\qquad\qquad\qquad\qquad\qquad
$\boxed{\boxed{\aligned
&\text{\eightpoint"Thalassa Aitheria" Reports}\\
&\text{\eightpoint RCMPI-95/01}\endaligned}}$\newline
\ \newline
\ \newline
\ \newline
\topmatter
\title On the nonHamiltonian interaction of two rotators
\endtitle
\author Denis V. Juriev
\endauthor
\affil Research Center for Mathematical Physics\linebreak
and Informatics "Thalassa Aitheria"\linebreak
ul.Miklukho-Maklaya 20-180, Moscow 117437 Russia\linebreak
(e mail: denis\@juriev.msk.ru)\linebreak
\endaffil
\abstract Classical dynamical equations describing a certain version of the
nonHamiltonian interaction of two rotators (Euler tops with completely
degenerate inertia tensors) are considered. The simplest case is integrated.
It is shown that the dynamics is almost periodic with periods depending on
the initial data.
\endabstract
\endtopmatter

Classical dynamics of Hamiltonian systems is often described by remarkable
algebraic structures such as Lie algebras [1] or their nonlinear
generalizations [2]. There is a hope that not less important algebraic
objects govern a behaviour of the interacting Hamiltonian systems. It seems
that these mathematical objects may be unravelled in a certain formal way.

There exist several types of an interaction of Hamiltonian systems:
often it has a potential character, sometimes it is ruled by a deformation of
the Poisson brackets; however, one of the most intriguing and less explored
forms is a nonHamiltonian interaction. In general, it can not be described
by deformations of the standard Hamiltonian data (Poisson brackets and
Hamiltonians). Sometimes, such interaction is realized by the dependence of
the Poisson brackets of one Hamiltonian system on the state of another.
Dynamical systems of charged objects with nonpotential magnetic interactions
provide a realistic and relatively simple example (cf.[3]). More sophisticated
generalizations may be also formulated (e.g. two nucleons at the classical
level).

To unravel certain mathematical structures beyond such interactions it is
reasonable to consider some simple toy constructions. The examples of a linear
dependence of Poisson brackets on the states were considered in [4]. Physically
realistic example is a pair of spherical charged particles with fixed centers
interacting via Biot--Savart--Amp\'ere forces. However, if charged particles
are asymmetric then one should consider each particle in its proper
coordinates.
In such a way equations of motion of each object will contain a nonlinear
dependence on the states of an opposite one.

A toy example of the quadratic nonlinear dependence will be considered below.
Some physical systems with magnetic interactions are characterized by such
dependence, e.g. two pairs of charged spherical particles. In this case the
dependence is linear as in velocities as in distances between particles in
the pairs. Some features of our toy example (boundness of motion, stability
and almost periodicity) may be also peculiar to a class of algebraically
constructed systems with quadratic nonHamiltonian magnetic--type interactions.

Our general constructions may be useful for a qualitative analysis of
essentially nonlinear systems (cf. [2]), because in many cases the quadratic
approximation determines the qualitative behavior of interacting systems, so
the highest perturbations do not disturb the boundness and stability of motion.

First, let us specify the necessary algebraic requisites.

\definition{Definition 1} A pair $(V_1,V_2)$ of linear spaces will be called
{\it I--pair\/} iff there are defined (nonlinear) mappings
$$\aligned &h_1:V_2\mapsto\Hom(\Lambda^2(V_1),V_1)\\
&h_2:V_1\mapsto\Hom(\Lambda^2(V_2),V_2)
\endaligned$$
such that $\forall A\in V_2$ its image $h_1(A)\in\Hom(\Lambda^2(V_1),V_1)$ is
a Lie bracket in $V_1$ and $\forall X\in V_1$ its image
$h_2(X)\in\Hom(\Lambda^2(V_2),V_2)$ is a Lie bracket in $V_2$. The Lie
bracket in $V_2$ corresponded to $A$ will be denoted by $[\cdot,\cdot]_A$,
whereas the Lie bracket in $V_1$ corresponded to $X$ will be denoted by
$[\cdot,\cdot]_X$.
\enddefinition

Here $\Hom(H_1,H_2)$ denotes the space of all linear operators from $H_1$ to
$H_2$, $\Lambda^2(H)$ is a skew square of the linear space $H$, so
$\Hom(\Lambda^2(H),H)$ is the space of all skew--symmetric bilinear binary
operations in $H$. The Lie brackets in $H$ form a submanifold
$\operatorname{Lie}(H)$ of the space $\Hom(\Lambda^2(H),H)$.

Mechanically, $V_1$ and $V_2$ are (possibly somehow reduced) phase spaces of
the interacting systems.

\

The examples of I--pairs with linear mappings $h_i$ were considered in [4].

\remark{Example of a nonlinear I--pair} Let $V_i$ be the adjoint
representations
of the Lie algebra $\so(n)$: $V_i\simeq\Lambda^2(\Bbb R^n)\subseteq\Mat_n(\Bbb
R)$. Then put
$$[X,Y]_A=XY-YX+\varepsilon(XA^2Y-YA^2X)$$
for all $X,Y\in V_1$ and $A\in V_2$ as well as
$$[A,B]_X=AB-BA+\varepsilon(AX^2B-BX^2A)$$
for all $A,B\in V_2$ and $X\in V_1$.
\endremark

\definition{Definition 2} Let us consider two elements $\Omega^+$ and
$\Omega^-$
in $V_1$ and $V_2$, respectively. The equations
$$\dot X_t=[\Omega^+,X_t]_{A_t},\qquad \dot A_t=[\Omega^-,A_t]_{X_t},$$
where $X_t\in V_1$ and $A_t\in V_2$ are called {\it the (nonlinear) dynamical
equations\/} ({\it Euler formulas\/}) {\it associated with the I-- pair
$(V_1,V_2)$ and the elements $\Omega^{\pm}$}.
\enddefinition

Such formulas generalise classical Euler formulas for rotator [5, p.129,136].

The toy examples for the linear case was considered in [4]. Physically
realistic
examples of the linear dependence include magnetic interaction between
two charged symmetric Euler tops ($I_1=I_2$, $I_3=0$, $I_n$ are components of
inertia tensor), two Cooper pairs or two positroniums.

The dynamical equations associated with the I--pair of the example above and
with $\Omega^+\in V_2$ and $\Omega^-\in V_1$ may be considered as describing
the nonHamiltonian interaction of two n--dimensional rotators: $\varepsilon$
is the coupling constant.

\remark{Lemma} If $\Omega^{\pm}=\Omega$ then the system has two independent
integrals of motion: $\Cal I=\left<A,A\right>+\left<B,B\right>$ and
$\Cal J=\left<\Omega,A\right>^2+\left<\Omega,B\right>^2$ ($\left<X,Y\right>=
\operatorname{Tr}(XY)$).
\endremark

\

Really, $\Cal I$ should be identified with energy so the lemma means that
the system is not dissipative.

\

Let's consider the case $n=3$. In this case the commutator $[\cdot,\cdot]$ in
$\so(3)$ reduces to an ordinary vector product $\times$ in $\Bbb R^3$.

The dynamical equations for two independent rotators have the form
$$\left\{\aligned
\dot A&=[\Omega^+,A]\\
\dot B&=[\Omega^-,B]
\endaligned\right.$$
$A, B, \Omega^{\pm}\in\so(3)$ or
$$\left\{\aligned
\dot A&=\Omega^+\times A\\
\dot B&=\Omega^-\times B
\endaligned\right.$$
$A, B, \Omega^{\pm}\in\Bbb R^3$.

If $A$ and $B$ are identified with coordinates they describe the motion of
a rigid body with fixed axis. If $A$ and $B$ are identified with velocities
they describe a charged particle in an external magnetic field or a rotator
[5, p.129,136]. If $A$ and $B$ are identified with angular velocities they
describe a symmetric Euler top [5, p.128, 136, 142, 146] (see also
[6,Sect.I,Ch.1]).

The perturbed equations have the form
$$\left\{\aligned
\dot A&=[\Omega^+,A]+\varepsilon(\Omega^+B^2A-AB^2\Omega^+)\\
\dot B&=[\Omega^-,B]+\varepsilon(\Omega^-A^2B-BA^2\Omega^-)
\endaligned\right.$$
Here $\so(3)$ is canonically imbed into $\Mat_n(\Bbb R)$ as a space of
skew--symmetric matrices.

\proclaim{Theorem} The dynamical equations for two interacting rotators are
integrable for n=3 and $\Omega^+=\Omega^-=\Omega$. A motion defined by them
is almost periodic with periods depending on the initial data.
\endproclaim

\demo{Proof}
Let $\bold i$, $\bold j$, $\bold k$ be a canonical orthonormal basis in
$\so(3)\simeq\Bbb R^3$:
$[\bold i,\bold j]=\bold i\times\bold j=\bold k$,
$[\bold j,\bold k]=\bold j\times\bold k=\bold i$,
$[\bold k,\bold i]=\bold k\times\bold i=\bold j$.
Put $\Omega=\Omega^{\pm}=\omega\bold i$ and expand $A$ and $B$ as
$a_i\bold i+a_j\bold j+a_k\bold k$ and $b_i\bold i+b_j\bold j+b_k\bold k$,
then the dynamical equations will be written as
$$\left\{\aligned
\dot a_i&=\varepsilon\omega b_i(b_ja_k-a_jb_k)\\
\dot b_i&=\varepsilon\omega a_i(a_jb_k-a_kb_j)
\endaligned\right.$$

$$\left\{\aligned
\dot a_j&=-\omega a_k+\varepsilon\omega b_j(b_ja_k-a_jb_k)\\
\dot b_j&=-\omega b_k+\varepsilon\omega a_j(a_jb_k-a_kb_j)\\
\dot a_k&=\omega a_j+\varepsilon\omega b_k(b_ja_k-a_jb_k)\\
\dot b_k&=\omega b_j+\varepsilon\omega a_k(a_jb_k-a_kb_j)
\endaligned\right.$$

Put
$$\aligned a_j=a_0\cos\varphi,& a_k=a_0\sin\varphi\\
b_j=b_0\cos\psi,& b_k=b_0\sin\psi
\endaligned$$
then
$$\left\{\aligned
\dot a_i&=\varepsilon\omega b_ia_0b_0\sin(\varphi-\psi)\\
\dot b_i&=-\varepsilon\omega a_ia_0b_0\sin(\varphi-\psi)
\endaligned\right.$$

Put $a_i^2+b_i^2=\Cal A^2$, $a_i=\Cal A\cos\chi$, $b_i=\Cal A\sin\chi$, then
$$\dot\chi=-\varepsilon\omega a_0b_0\sin(\varphi-\psi).$$

Note that
$$\left\{\aligned
(a_0^2)^\cdot&=\varepsilon\omega a_0^2b_0^2\sin 2(\varphi-\psi)\\
(b_0^2)^\cdot&=-\varepsilon\omega a_0^2b_0^2\sin 2(\varphi-\psi)
\endaligned\right.$$

Put $a_0^2+b_0^2=\Cal B^2$, $a_0=\Cal B\cos\xi$, $b_0=\Cal B\sin\xi$, then
$$\left\{\aligned
\dot\xi&=-\tfrac{\varepsilon\omega\Cal B^2}4\sin 2(\varphi-\psi)\sin 2\xi\\
\dot\chi&=-\tfrac{\varepsilon\omega\Cal  B^2}2\sin
2\xi\sin\eta\endaligned\right.$$

Moreover,
$$\left\{\aligned
\dot\varphi&=\omega-\varepsilon\omega b^2_0\sin^2 2(\varphi-\psi)\\
\dot\psi&=\omega-\varepsilon\omega a^2_0\sin^2 2(\varphi-\psi)
\endaligned\right.$$

Denote $\vartheta=\tfrac12(\varphi+\psi)$, $\eta=\varphi-\psi$. Then
$$\left\{\aligned
\dot\vartheta&=\omega+\tfrac{\varepsilon\omega\Cal B^2}2\sin 2\eta\\
\dot\eta&=-\varepsilon\omega\Cal B^2\sin^2\eta\cos 2\xi
\endaligned\right.$$

Mark that $K=\sin\eta\sin 2\xi$ is an integral of motion, therefore,
$$\left\{\aligned\dot\eta=&\varepsilon\omega\Cal B^2\sin^2\eta\cos 2\xi\\
\dot\xi=&-\tfrac{\varepsilon\omega\Cal B^2}4\sin 2\eta\sin 2\xi
\endaligned\right.$$
whereas
$$\dot\theta=\omega-\tfrac{\varepsilon\omega\Cal B^2}2\sin^2\eta.$$

The system of dynamical equations on $\eta$ and $\xi$ is easily integrated.
Namely,
$$\cos 2\xi=\sqrt{1-K^2}\cos(C_1+\varepsilon\omega\Cal B^2Kt)$$
and
$$\cot\eta=\tfrac{\sqrt{1-K^2}}K\sin(C_1+\varepsilon\omega\Cal B^2Kt)).$$

Hence,
$$\theta=\omega t-\tfrac{\varepsilon\omega\Cal B^2K^2}2\int
\frac{dt}{1-(1-K^2)\cos^2(C_1+\varepsilon\omega\Cal B^2Kt)}+C_2.$$
The last integral may be easily computed but the resulted expression is not
much useful, one should only mark that $\theta=\omega't+f(\varepsilon
\omega\Cal B^2Kt)$, where $f$ is a $2\pi$--periodic function, $\omega'=
\omega+\varepsilon\omega\Cal B^2K\Delta$ ($\Delta=\Delta(K)$).
One should mark also that
$$\chi=-\tfrac{\varepsilon\omega\Cal B^2K}2t+C_3,$$
so the periods of the almost periodic motion are $\tfrac{2\pi}{\omega'}$ and
$\tfrac{4\pi}{\varepsilon\omega\Cal B^2K}$.
\qed\enddemo

Note that the systems admits three integrals of motion: $\Cal I$, $\Cal J$
and $K$. The presence of two of them is provided by the algebraic structure
of the model (and equality $\Omega^+=\Omega^-$; see above). It seems that the
third integral has no any natural algebraic interpretation.

\

Thus, the classical dynamical equations describing a certain toy version of the
nonHamiltonian interaction of two rotators are considered. The simplest case
of the dimension $n=3$ and the coinciding angular velocities
$\Omega^{\pm}=\Omega$ is integrated. It is shown that a motion is almost
periodic with periods depending on the initial data.

One may suppose that a qualitative behavior of the considered toy system (e.g.
boundness and almost periodicity) is peculiar to a more general class of
Hamiltonian systems with nonHamiltonian interaction of magnetic type. Thus, the
results may be of a possible use for a qualitative and structural analysis of
complicated physical systems of such kind.

However, note that a relation between the algebraic construction of the model
and its explicit integration is formally almost "accidental". Nevertheless,
it seems that this relation is essential. It is reasonable to formulate a
{\bf hypothesis} that if two Hamiltonian systems are integrable and their
integrability is provided by some algebraic mechanism, then an integrability of
the interacting dynamics may be related to some additional new algebraic
structures, which underlie the interaction.

Certainly, such "phenomenological" hypothesis is merely a constructive
compromise between practical needs and theoretical abilities. It does not
explain a mechanism of the relation but only suggest to look for any explicit
integration of algebraically constructed nonHamiltonian systems. In such a way
it may be useful.

\newpage

\Refs
\roster
\item"[1]" Arnold V I 1976 {\it Mathematical methods of classical mechanics},
Springer--Verlag; Dubrovin B A, Novikov S P and Fomenko A T 1988 {\it Modern
geometry}, Springer--Verlag; Fomenko A T 1988 {\it Symplectic geometry.
Methods and applications}, Gordon and Beach, New York.
\item"[2]" Karasev M V and Maslov V P 1993 {\it Nonlinear Poisson brackets.
Geometry and quantization}, Amer.~Math.~Soc., RI.
\item"[3]" Dubrovin B A, Krichever I M and Novikov S P 1988 {\it Current
Problems Math, Fundamental Directions IV}, Moscow, VINITI.
\item"[4]" Juriev D 1994 {\it Russian J Math Phys\/} 3(4) [e--print version
(LANL Electronic Archive on Solv. Integr. Systems): {\it solv-int/9505001}];
Classical and quantum dynamics of noncanonically coupled oscillators and Lie
superalgebras. E--print (SISSA Electronic Archive on Funct. Anal.): {\it
funct-an/9409003}.
\item"[5]" Landau L D and Lifschitz E M 1958 {\it Mechanics}. Fizmatlit,
Moscow.
\item"[6]" Leimanis E 1958 {\it Surveys Appl Math\/} 2. John Wiley \& Sons
Inc.,
Ney York; Chapman \& Hall Ltd., London.
\endroster
\endRefs
\enddocument